\documentclass[11pt,a4paper]{amsart} 
\usepackage[english]{babel}
\usepackage{amssymb,xspace,enumerate,graphics,ifthen}
\usepackage{amstext}
\theoremstyle{plain}
\usepackage{amsbsy,amssymb,amsfonts,latexsym}

\date{\today}
\title{$J$-holomorphic discs and real analytic hypersurfaces}
\author{William ALEXANDRE}
\address{Laboratoire Paul Painlev\'e U.M.R. CNRS 8524, U.F.R. de
Math\'ematiques,  cit\'e scientifique, Universit\'e Lille 1, F59 655 Villeneuve d'Ascq Cedex, France.}
\email{ william.alexandre@math.univ-lille1.fr}
\thanks{The first author is partially supported by A.N.R. BL-INTER09-CRARTIN}
\author{Emmanuel MAZZILLI}
\address{Laboratoire Paul Painlev\'e U.M.R. CNRS 8524, U.F.R. de
Math\'ematiques,  cit\'e scientifique, Universit\'e Lille 1, F59 655 Villeneuve d'Ascq Cedex, France.}
\email{ emmanuel.mazzilli@math.univ-lille1.fr}
\date{}
\newtheorem{theorem}{Theorem}[section]
\newtheorem{lemma}[theorem]{Lemma}

\newtheorem{example}[theorem]{Example}
\newtheorem{proposition}[theorem]{Proposition}
\newtheorem{corollary}[theorem]{Corollary}
\newtheorem{remark}{\it Remark}
\newtheorem{definition}[theorem]{Definition}
\newtheorem{notation}[theorem]{Notation}
\def \pint {\vbox{ \hbox to 5 pt {\hfil \vrule height 4pt}\hrule}\hskip 3pt}
\def\leqs{\lesssim}

\def\cc{\mathbb{C}}
\def\ll{\mathbb{L}}
\def\rr{\mathbb{R}}
\def\nn{\mathbb{N}}
\def\dd{\mathbb{D}}

\def\re{{\rm Re} }
\def\im{{\rm Im} }
\def \qed {\hbox{\hskip 5pt} \vbox{\hrule \hbox to 5pt 
{\vrule height 4.2pt \hfil \vrule}\hrule}}
\def \pint {\vbox{ \hbox to 5 pt {\hfil \vrule height 4pt}\hrule}\hskip 3pt}

\newcommand{\cal}{\mathcal}

\newcommand{\app}[5]
{#1: \left \{ 
\begin {array}{ccl} 
#2 & \longrightarrow & #3 \\#4 & \longmapsto & #5 \end {array} \right .}

\newcommand{\diffp}[2]{\frac{\partial #1}{\partial #2}}
\newcommand{\mlabel}[1]{\label {#1}}
\renewcommand{\over}[2]{\genfrac{}{}{0pt}{}{#1}{#2}}

\newcommand{\pk}[2]{{\cal P}_{{\ifthenelse{#1=1}{ }{#1}}\kappa|\rho({#2})|}({#2})}

\newcommand{\pr} {\noindent{\it Proof:} }

\renewcommand{\aa}{{\cal A}}
\newcommand{\cll}{{\cal L}}
\newcommand{\x}[1]{{\frac{\partial}{\partial x_{#1}}}}
\newcommand{\y}[1]{{\frac{\partial}{\partial y_{#1}}}}
\def \indic{\mathbf{1}}
\begin{document}
\pagestyle{plain}

\begin{abstract}
We give in $\rr^6$ a real analytic almost complex structure $J$, a real analytic hypersurface $M$ and a vector $v$ in the Levi null set at $0$ of $M$, such that there is no germ of $J$-holomorphic disc $\gamma$ included in $M$ with $\gamma(0)=0$ and $\diffp\gamma x(0)=v$, although the Levi form of $M$ has constant rank. Then for any hypersurface $M$ and any complex structure $J$, we give necessary conditions under which there exists such a germ of disc.
\end{abstract}
\maketitle
\section{Introduction}
Throughout all this paper, we denote by $M$ be a real analytic hypersurface in $\rr^{2n}$ such that $0$ belongs to $M$. We denote by $\varphi $ a real analytic function such that $M=\{z\in\rr^{2n},\ \varphi(z)=0\}$ and such that $d\varphi$ does not vanish on $M$.  We also denote by 
$TM$ the tangent vector bundle to $M$.\\
We equip $\rr^{2n}$ with a real analytic almost complex structure $J$, i.e.  a linear map from $T\rr^{2n}$ into itself such that $J^2=-Id$. Without restriction, we can always assume that $J(0)=J_0$, the standard complex structure in $\rr^{2n}=\cc^n$, and that $J$ is as close as we need to $J_0$. We denote by $T^JM$ the $J$ invariant part of $TM$, that is $T^JM=TM\cap J\ TM$.

The Levi form $\cll_\varphi $ of $M$ and its kernel are two objects linked to the geometry of $M$. The Levi form is defined for every vector field $X\in T^JM$ by 
$${\cal L}_\varphi(X)= d\varphi(J[X,JX]).$$
We will also consider $\hat{\cal L}_\varphi$, the polar form of $\cll_\varphi$:
$$\hat\cll_\varphi(X,Y)=d\varphi(J[X,JY]+J[Y,JX])+id\varphi(J[X,Y]+J[JX,JY]).$$
Both $\cll_\varphi $ and $\hat\cll_\varphi $ do not depend on the derivatives of $X$ and $Y$. Moreover, they depend on the defining function $\varphi$ up to multiplication by a nonnegative function. Therefore $\ker\cll_\varphi $ does no depend on $\varphi $ and we simply denote it by $\ker\cll$.
\par\medskip
When $J$ is integrable, or more simply when $J$ is the standard structure, and when $\ker\cll$ is a subbundle of $T^JM$, Freeman proved in \cite{Fre1} that  there exists a complex foliation of $M$. However, when $J$ is a generic almost complex structure, there is no foliation by $J$-holomorphic manifolds of dimension greater than or equal to 2 because such objects do not exist in general. Hence Freeman's theorem does not hold anymore in the non integrable case; we can only hope to find a $J$-holomorphic disc with prescribed derivative of order 1 included in $M$. A $J$-holomorphic disc is a $C^1$ map $\gamma:\dd\to\rr^{2n}$, where $\dd$ is the unit disc of $\cc$ equipped with the standard structure $J_0$,  such that for all $\zeta\in\dd$, we have as linear map:
$$d\gamma (\zeta) \circ J_0=J(\gamma (\zeta)) \circ d\gamma (\zeta).$$

When $\gamma$ is a $J$-holomorphic disc included in $M$, i.e.\ $\gamma(\dd)\subset M$, the vector  $\diffp{\gamma}{x}(0)$ belongs to $(\ker \cll)_{\gamma(0)}$. When $J$ is integrable and $\ker \cll$ is a subbundle of $T^JM$, Freeman's result implies that the converse is true (even in the non real analytic case); that is, for any point $p\in M$ and any $v\in\ker\cll$, there exists a $J$-holomorphic disc $\gamma$ included in $M$ such that $\gamma(0)=p$ and $\diffp\gamma x(0)=v$. This result was generalized by Kruzhilin and Sukhov in \cite{KS} when $J$ is not integrable and  $M$ is Levi flat, that is when $\ker\cll=T^JM$. It was then natural to hope that this result was true in the real analytic case when $\ker\cll$ is only a subbundle of $T^JM$. However, we will show that this is not the case, even when $M$ is pseudoconvex. We will prove in Section \ref{contre-exemple} the following theorem:
\begin{theorem}\label{th0}
There exist a real analytic almost complex structure $J$, a real analytic hypersurface $M$ in $\rr^6$ such that $\ker\cll$ is a subbundle of dimension $1$ of $T^JM$ and  a vector $v\in \ker\cll_0$ such that no $J$-holomorphic disc $\gamma$, satisfying  $\diffp\gamma  x(0)=v$, is tangent to $M$ at $0$ at order greater than $5$.
\end{theorem}

The following natural question then arises: under which condition on $v\in T^J_0M$ does there exist a $J$-holomorphic disc $\gamma$ included in $M$ such that $\gamma(0)=0$ and $\diffp{\gamma}{x}(0)=v$ ?

In order to answer this question, we fix a symmetric connexion $\nabla$ and for simplicity's sake, we use a multiplicative notation: when $X$ and $Y$ are two vector fields,  we denote by $X.Y$ the vector field $\nabla_X Y$. We also introduce the following notations:
\begin{notation}
If $X$ is a vector field and $p$ and $q$ are two integers, we define the vector field $D_X^{p,q}$ by the formula
$$D_X^{k,l}\phi=\underbrace{JX\cdot(JX\cdots (JX}_l\cdot \underbrace{(X\cdot(X\ldots (X}_k\cdot\phi)\ldots).$$
\end{notation}
\begin{notation}
 Let $X$ be a vector field. We denote by $\ll(X)$ the complex Lie algebra generated by $X$, i.e.\ the smallest complex Lie algebra which contains $X$ and $JX$ and is such that for $Y$ and $Z$ in $\ll(X)$, $[Y,Z]$ belongs to $\ll(X)$.
\end{notation}

We will establish two sufficient conditions on a vector $X$ under which there exists a $J$-holomorphic disc $\gamma$ included in $M$ such that for all $k\in\nn^*$, $\diffp{^k\gamma}{x^k}(0)=D_X^{k-1,0} X(0)$. The first one is given by the following theorem:
\begin{theorem}\label{alg-lie}
Let $X\in T^JM$ be a real analytic vector field in a neighborhood of $0\in M$ such that $\ll(X)$ is included in $T^JM$.\\
Then, for all $p\in M$ in a neighborhood of $0$, there exists a germ of $J$-holomorphic disc $\gamma$ such that 
\begin{enumerate}
 \item $\gamma(0)=p$,
 \item for all $k\in\nn^*$, $\diffp{^k\gamma}{x^k}(0)=D_X^{k-1,0} X(p)$.
\end{enumerate}
\end{theorem}
Unlike Freeman's Theorem \cite{Fre1}, Theorem \ref{alg-lie} makes no assumption on the rank of the Levi form $\cll$. In particular, Theorem \ref{alg-lie} can be applied to an integrable complex structure when $\ker\cll$ is not a subbundle. In this case, we will show that for any $X\in\ker\cll$, there exists a $J$-holomorphic disc $\gamma $ included in $M$ such that $\diffp\gamma x(0)=X(0)$  (see Theorem \ref{freeman}). We will also give an example of application of our result in a non-integrable case and where the Levi form has constant rank (see Example \ref{ex1}).

In order to state our next result which requires a condition in the spirit of Freeman's theorem, we recall the following definition introduced in \cite{BM}:
\begin{definition}\label{def-commute}
We say that a vector field $X$ commutes at order $k$ at a point $p$ if for all $l\leq k$, all $X_1,\ldots, X_l\in\{X,JX\}$, 
$$[X_1,[\ldots, [X_{l-1},X_l]\ldots](p)=0.$$
\end{definition}
We will prove the following result:
\begin{theorem}\label{subbundle}
Let $\mathbb{L}$ be a subbundle of $T^JM$ such that the following property is true: if $X\in \ll$ commutes at order $k$ at the point $0$, then for all $X_1,\ldots, X_{k+1}\in \{X,JX\}$, $[X_1,[\ldots,[X_{k},X_{k+1}]\ldots]](0)$ belongs to $\ll_0$.\\
Then for all $X\in \ll$, there exists a germ of $J$-holomorphic disc $\gamma $ such that $\gamma(0)=0$, $\diffp{^k\gamma }{x^k}(0)=D_X^{k-1,0}X(0)$ for all $k\in\nn^*$ and $\gamma (\dd)\subset M$. 
\end{theorem}
This Theorem \ref{subbundle} enables us to prove, in the real analytic case and in a constructive way, Kruzhilin and Suhkov's theorem \cite{KS} (See Theorem \ref{KS}). We will also give examples of applications of Theorem \ref{subbundle} (see Example \ref{ex2}).
\par\medskip
In order to prove Theorem \ref{alg-lie} and \ref{subbundle}, we prove a generalization of the existence results of $J$-holomorphic discs with prescribed jets of finite order, due to Sikorav in the case of $1$-jets (see \cite{Sik}) and to Ivashkovich and Rosay \cite{IR} and independently to Barraud and Mazzilli \cite{BM} in the case of $k$-jets, for $k$ finite.
For this, a key point is to show that if $(x_k)_{k\in\nn}$ is an infinite sequence of vectors in $\rr^{2n}$ such that the series $\sum_{k} \frac{k+1}{k!}|x_k|$ converges, then there is a germ of $J$-holomorphic disc $\gamma$ included in $M$ such that $\diffp{^k\gamma}{x^k}(0)=x_k$ for all $k\in\nn$ (see Theorem \ref{jhpd} and Corollary \ref{gjhpd}). The general procedure to prove this is in some sense analogue to the case of finite jets but the proof is more intricate: we have to prove that some maps between two Banach spaces $\tilde\Omega$ and $\Omega$ are invertible in a neighborhood of the origin. In the case of $k$-jets, the Banach spaces considered were simply $\tilde\Omega=(\rr^{2n})^k$ and $\Omega=C^{k+\alpha}(\overline{\dd})$ for $\alpha>0$, 
with the $k$-jet belonging to $(\rr^{2n})^k$ and the disc to $C^{k+\alpha}(\overline{\dd})$. 
The first difficulty we have to overcome is to find the appropriate Banach spaces.
Here,  the $k$-jet becomes an infinite sequence of vectors of $\rr^{2n}$, so we have to consider a space of sequences with an appropriate norm. This role will be played by the space $\tilde\Omega$  of all sequences $(x_k)_k$ of vectors in $\rr^{2n}$ such that the norm $\|(x_k)_k\|_{\sim{}}:= \sum_{k} \frac{k+1}{k!}|x_k|$ is finite. The space of functions will be the space $\Omega$ of real analytic functions $f=\sum_{k,l} f_{k,l}z^k\overline z^l$ on the disc $\dd$  such that $\|f\|:=\sum_{k,l}(1+k+l+kl)|f_{kl}|$ is finite. Then, as in the case of $k$-jets, we associate to any jet a function and to any function a jet (See Theorem \ref{jhpd} for details). Showing that these maps are firstly well defined and secondly invertible, we find a $J$-holomorphic disc with prescribed derivatives of any order.

In Proposition \ref{evf} we will prove that if $X$ is a real analytic vector field in a neighborhood of $0$, then the sequence $(D_X^{k-1,0}X(0))_{k\in\nn^*}$ belongs to $\tilde\Omega$. Hence there exists a germ of $J$-holomorphic disc $\gamma$ such that $\gamma(0)=0$ and $\diffp{^k\gamma}{x^k}(0)=D_X^{k-1,0}X(0)$ for all $k\in\nn^*$.
\par\medskip
From this point, the problem of finding a $J$-holomorphic disc included in $M$ with prescribed derivatives is reduced to finding a sequence of vector fields $(X_l)_l$ in $T^JM$ such that the sequence $(D^{k-1,0}_{X_l} X_l(0))_{k\in\nn^*}$ does not depend on $l$ and is such that $X_l$ commutes at order $l$ at $0$. Indeed, Theorem \ref{jhpd} gives a $J$-holomorphic disc $\gamma$ such that for all $k$ and all $l$, $\diffp{^k\gamma} {x^k}(0)=D^{k-1,0}_{X_l} X_l(0)$. Since $X_l$ commutes at order $l$, Theorem 1 of \cite{BM} implies that $\gamma$ is tangent to $M$ at $0$ at order $l+1$, for all $l$. As both $M$ and $\gamma$ are real analytic, it follows that $\gamma$ is actually included in $M$. 

Starting from a vector field $X$ which does not commute but which satisfies the assumptions of Theorem \ref{alg-lie} or \ref{subbundle}, we inductively construct  such a sequence of vector fields by solving systems of linear equations (see Lemma \ref{pfreeman} and Lemma \ref{ppfreeman}).
\par\medskip
The paper is organized as follows. In section \ref{contre-exemple}, we prove Theorem \ref{th0}. In Section \ref{secjhpd}, we show the existence of $J$-holomorphic disc with prescribed derivatives at any order. Finally, in Section \ref{secjhhs}, we prove Theorem \ref{alg-lie} and \ref{subbundle} and  give examples of applications of these theorems.

\section{A counter-example to a generalization of Freeman's Theorem}\label{contre-exemple}
In this section, in order to prove Theorem \ref{th0},  we exhibit a complex structure $J$ and a pseudoconvex hypersurface $M$ in $\rr^6$, both real analytic, such that there exists a vector  $v$ which belongs to the kernel of the Levi form of $M$ at the point $0$ but which is not the derivative at $0$ of some $J$-holomorphic disc $\gamma:\dd\to \rr^{2n}$ included in $M$, despite the fact that the kernel of the Levi form of $M$ is a subbundle of $T^JM$. We will prove that in fact, there exists no $J$-holomorphic disc $\gamma$ tangent to $M$ at $0$ at order $5$ (i.e.\ such that $\varphi\circ\gamma(\zeta)=O(\zeta^5)$).
\par\medskip
{\it Proof of Theorem \ref{th0}:} 
Let $\varphi:\rr^6\to\rr$ be the map defined by $\varphi(x_1,y_1,x_2,y_2,x_3,y_3)=y_1$, and let $M$ be the set $M=\{z\in\rr^6\:\ \varphi(z)=0\}$.\\
We define the six following vector fields 
\begin{align*}
 L_1&=\diffp{}{x_1},	&L_2&=\y{1},\\
L_3&=\x2-\frac12y^2_3\diffp{}{x_1},&L_4&=\diffp{}{y_2}+(-2y_3x_3+x_2)\diffp{}{x_1},\\
L_5&=\x{3}-y_2y_3\x{1},		&L_6&=\y3+x_3\x2-\left(\frac{x_3y_3^2}2+x_3y_2\right)\x1,
\end{align*}
and the complex structure $J$ they induce by setting 
\begin{align*}
 JL_1&=L_2,&JL_3&=L_4,&JL_5&=L_6.
\end{align*}
Therefore $J(0)$ is the standard complex structure $J_0$ and $TM$ is spanned over $\rr$ by $L_1,$ $L_3,$ $L_4,$ $L_5$ and $L_6$, while $T^JM$ is spanned over $\cc$ by $L_3$ and $L_5$.\\
When we compute the Lie brackets of the $L_i$'s belonging to the complex tangent bundle, we get:
\begin{align*}
 [L_3,L_4]&=L_1,\\	
 [L_3,L_5]&=0,		&[L_4,L_5]&=y_3\x1,\\
 [L_3,L_6]&=y_3\x1	& [L_4,L_6]&=0, 	& [L_5,L_6]&=L_3.\\
\end{align*}
Therefore
\begin{align*}
 \hat \cll_\varphi(L_3,L_5)&=0,\\
 \hat\cll_\varphi(L_3,JL_3)&= d\varphi(L_2)=1,\\
 \hat \cll_\varphi(L_5,JL_5)&=d\varphi(L_3)=0,\\
 \end{align*}
which thus implies  that $M$ is pseudoconvex and that $\ker \cll=span_\cc\{L_5\}$; therefore $\ker \cll$ is a subbundle of $T^JM$.

Now let us assume that there exists a $J$-holomorphic disc $\gamma:\dd\to\rr^{2n}$ tangent to $M$ at order $5$ and such that $\diffp{\gamma}{x}(0)=L_5(0)$. Theorem 1 from \cite{BM} implies that there exists a vector field $X\in T^JM$ such that $X$ commutes at order 4 at $0$, i.e. such that all the Lie brackets of $X$ and $JX$ of length at most 4 vanish at $0$, and such that $X(0)=\diffp\gamma x(0)=L_5(0)$. We prove that such a vector field does not exist.\\
Let us consider $X=aL_3+bL_4+cL_5+dL_6$ where $a$, $b$, $c$ and $d$ are real valued functions. Let us assume that $X$ commutes at order 4 and that $X(0)=L_5(0)$. Therefore we have $a(0)=b(0)=d(0)=0$ and $c(0)=1$. We compute $[X,JX]$:
\begin{align*}
 [X,JX]=&(a^2+b^2) L_1+(c^2+d^2)L_3+(-X(b)-JX(a)) L_3 + (X(a)-JX(b))L_4\\&+(-X(d)-JX(c))L_5+(X(c)-JX(d))L_6
\end{align*}
and since $[X,JX](0)=0$, it follows that
\begin{align}
    X(b)(0)+JX(a)(0)&=-1.\mlabel{eq1}
\end{align}
We also get from the computation of $[X,JX]$ that
$$d\varphi(J[X,JX])=a^2+b^2,$$
thus
$$X\cdot X\cdot d\varphi(J[X,JX])=2(X(a)^2+X(b)^2+aX\cdot X(a)+b X\cdot X(b))$$
and
$$JX\cdot JX\cdot d\varphi(J[X,JX])=2(JX(a)^2+JX(b)^2+aJX\cdot JX(a)+b JX\cdot JX(b)).$$
Using the commutativity of $X$ at $0$, we will show that $X\cdot X\cdot d\varphi(J[X,JX])(0)=
JX\cdot JX\cdot d\varphi(J[X,JX])(0)=0$ which is incompatible with (\ref{eq1}).

We have
\begin{align*}
 X\cdot d\varphi(J[X,JX])=&d\varphi((X.J)[X,JX])+d\varphi(J\cdot X\cdot [X,JX])+D^2\varphi(X,J[X,JX])
\end{align*}
and
\begin{align*}
&X\cdot X\cdot d\varphi(J[X,JX])\\
&=d\varphi((X\cdot X\cdot J)[X,JX])+
2d\varphi((X.J)\cdot X\cdot [X,JX])+
d\varphi(J(X\cdot X\cdot  [X,JX]))\\
&\hskip10pt+
D^2\varphi(X\cdot X,J[X,JX])+ D^2\varphi(X,(X\cdot J)[X,JX])+D^2\varphi(X,J(X\cdot [X,JX]))\\
&\hskip10pt+D^3\varphi(X,X,J[X,JX]).
\end{align*}
Now, $[X,JX]$ vanishes at $0$ so $[X,[X,JX]](0)=X\cdot [X,JX](0)$ and since $X$ commutes at $0$ at order 4, $X\cdot[X,JX](0)=0$.\\
Then we get
\begin{align*}
 [X,[X,[X,JX]]](0)&=X.[X,[X,JX]](0)\\
&=X\cdot( X\cdot [X,JX])(0)-X\cdot([X,JX]\cdot X)(0).
\end{align*}
Using again $[X,JX](0)=0$, we get $X\cdot([X,JX]\cdot X)(0)=(X\cdot [X,JX])\cdot X(0)$ and since $X\cdot [X,JX](0)=0$, we get $X\cdot ([X,JX]\cdot X)(0)=0$. Therefore $X\cdot(X\cdot[X,JX])(0)=[X,[X,[X,JX]]](0)$ and since $X$ commutes at order 4 at $0$, $X\cdot(X\cdot[X,JX])(0)=0$, which gives $X\cdot \left(X\cdot d\varphi( J[X,JX])\right)(0)=0$, and so $X(b)(0)=0$.

Analogously, we also have
$JX\cdot\left(JX\cdot d\varphi(J[X,JX])\right)(0)=0$ and so $JX(a)(0)=0$.

Now, just notice that (\ref{eq1}) and $X(b)(0)=JX(a)(0)=0$ are incompatible. This gives that there is no $J$-holomorphic disc $\gamma$ tangent at order $5$ at $0$ to $M$ such that $\diffp{\gamma}x(0)=L_5(0)$.

\section{Existence of $J$-holomorphic discs with prescribed derivatives}\label{secjhpd}
In this section, we show under an appropriate assumption, that there exists a $J$-ho\-lo\-morphic disc $\gamma$ with prescribed derivatives of any order. More precisely,  we prove the following
\begin{theorem}\mlabel{jhpd}
 Let $J$ be a real analytic complex structure in a neighborhood of the origin of $\rr^{2n}$, and let $(x_k)_k$ be a sequence of vectors of $\rr^{2n}$ such that $\sum_{k=0}^\infty \frac{k+1}{k!} |x_k|\leq r$, where $r>0$ is sufficiently small.\\
 Then there exists a $J$-holomorphic disc $\gamma:\dd\to\rr^{2n}$ such that for all $k\in\nn$, $\diffp{^k\gamma}{x^k}(0)=x_k$.
\end{theorem}
Before proving Theorem \ref{jhpd}, we rewrite the condition of $J$-holomorphicity of a disc when $J$ is close to the standard structure $J_0$. As in \cite{IS}, $\gamma:\dd\to\rr^{2n}$ is a $J$-holomorphic disc if and only if
$$\diffp{\gamma }{\overline\zeta}-A_J(\gamma )\diffp{\overline\gamma}{\overline\zeta}=0$$
with $A_J(z)=(J_0+J(z))^{-1}(J_0-J(z))C$, $C$ being the $\rr$-linear application which corresponds to the complex conjugation on $\rr^{2n}=\cc^n$.\\
Let $T:C^\omega(\dd)\to C^\omega(\dd)$ be the map defined by $T(u)(\zeta):=\int_{[0,\overline\zeta]} u(\zeta,\omega )d\omega.$ The real analytic function $T(u)$ is a primitive of $u$ with respect to $\overline{\zeta}$, that is $\diffp{T(u)}{\overline\zeta}=u$. Let also $\Phi_J:C^\omega(\dd)\to C^\omega(\dd)$ be the map defined by $\Phi_J(u)=u-T\left(A_J(u)\diffp{\overline u}{\overline\zeta}\right)$.\\
Then, a real analytic function $\gamma:\dd\to\rr^{2n}$ is $J$-holomorphic if and only if $\Phi_J(\gamma )$ is holomorphic in the classical way.
We can now prove Theorem \ref{jhpd}.

\noindent\it Proof of Theorem \ref{jhpd}: \rm The principle of the proof is analogous to the case of finite sequences: we first introduce two well chosen Banach algebras: an algebra $\tilde\Omega$ of vector sequences, and an algebra $\Omega$ of real analytic functions on $\dd$. We set 
\begin{align*}
 \tilde\Omega&=\left\{
  (x_k)_k\in\left(\rr^{2n}\right)^\nn,\ \left\|(x_k)_k\right\|_{\sim}:=
 \sum_{k=0}^\infty \frac{k+1}{k!}|x_k|<\infty\right\},\\
\Omega&=\left\{
f=\sum_{k,l}f_{k,l}\zeta^k\overline \zeta^l:\dd\to\rr^{2n},\ \|f\|:= \sum_{k,l=0}^\infty (1+k+l+kl)|f_{k,l}|<\infty\right\}.
\end{align*}
We also introduce two applications
\begin{align*}
& \app{\varphi_1}{\tilde\Omega}{C^\omega(\dd)}{(x_k)_k}{\sum_k \frac {x_k}{k!} \zeta^k},\\
&\app{\varphi_2}{\Omega}{\left(\rr^{2n}\right)^\nn}{f}{\left(\diffp{^kf}{x^k}(0)\right)_k}.
\end{align*}
We will prove the following facts:
\begin{enumerate}[(i)]
 \item\label{fact1} $\varphi _1$ is an isometry from $\tilde \Omega$  into $\Omega$,
 \item\label{fact2} $\varphi _2$ maps continuously $\Omega$ into $\tilde\Omega$,
 \item\label{fact3} $\Phi_J:\Omega\to\Omega$ is invertible in a neighborhood of  $0$,
 \item\label{fact4} $\varphi _2\circ \Phi_J^{-1}\circ \varphi _1:\tilde \Omega\to\Omega$ is also invertible in a neighborhood of $0$.
\end{enumerate}
Then, the disc $\gamma=\Phi_J^{-1}\circ \varphi _1\circ \left(\varphi _2\circ \Phi_J^{-1}\circ \varphi _1\right)^{-1}\left((x_k)_k\right)$ will be the  $J$-holomorphic disc we are looking for.\\
We first check fact (\ref{fact1}): $\varphi ((x_k)_k)=\sum_{k}\frac{x_k}{k!}\zeta^k$ and so, for all $l\neq 0$, the coefficient of $\zeta^k\overline \zeta^l$ in $\varphi _1((x_k)_k)$ vanishes. Therefore
$$\left\|\varphi_1((x_k)_k)\right\|=\sum_{k}(1+k)\frac{|x_k|}{k!}=\|(x_k)_k\|_\sim$$
and $\varphi _1$ is an isometry from $\tilde \Omega$  into $\Omega$.\\
We now check fact (\ref{fact2}): writing $\zeta$ as $\zeta=x+iy$, we get for all $k$ and $l$ that $\zeta^k\overline \zeta^l=x^{k+l}+yP_{k,l}(x,y)$ where $P_{k,l}$ is a polynomial in $x$ and $y$. Therefore
$$\left.\diffp{^n\zeta^k\overline \zeta^l}{x^n}\right|_{\zeta=0}=\begin{cases} (k+l)! &\text{ if } k+l=n\\0&\text{ otherwise}\end{cases}.$$
So, for $f=\sum_{k,l}f_{k,l}\zeta^k\overline\zeta^l$, we have
$\diffp{^nf}{x^n}(0)=n!\sum_{k+l=n}f_{k,l}$ and 
$$\varphi _2(f)=\left(n!\sum_{k+l=n} f_{k,l}\right)_{n\in\nn},$$
from which it follows that
\begin{align*}
 \|\varphi _2(f)\|_\sim&\leq \sum_n(n+1)\sum_{k+l=n}|f_{k,l}|\\
 &\leq \|f\|.
\end{align*}
So $\varphi _2(f)$ belongs to $\tilde \Omega$ and $\varphi _2:\Omega\to\tilde\Omega$ is continuous.\\
In order to prove Fact (\ref{fact3}), we introduce for $f\in\Omega$ the functions $g=A_J(f)\diffp{\overline f}{\overline z}$ and $G=T(g)$. Fact (\ref{fact3}) will be proved if we show that $G$ belongs to $\Omega$ and that $\|G\|<\frac12\|f\|$ when $\|f\|$ is small, because $\Phi_J$ will then be an invertible  perturbation of the identity in a neighborhood of $0\in \Omega$. We will use the following lemma.
\begin{lemma}\mlabel{crucial_lemma}
Let ${\cal A}$ be the Banach algebra of complex valued functions defined by
$${\cal A}=\{u(z,\overline z)=\sum_{k,l} u_{k,l}z^k\overline z^l,\ \|f\|_*:=\sum_{k,l}(1+k)|f_{k,l}|<\infty\},$$
and let $F$ be an analytic function bounded on $\{z\in\cc^k,\ |z|<R\}$, $R>0$. For $r>0$, we denote by ${\cal A}_r$ the set ${\cal A}_r=\{f\in {\cal A},\ \|f\|_*<r\}$.\\
Then, for all $r\in]0,R[$, the following properties hold true:
\begin{enumerate}[(i)]
 \item for all $f_1,\ldots, f_k\in {\cal A}_r,$ $F(f_1,f_2,\ldots, f_k)$ belongs to $\cal A$,
 \item $\app{\Psi_F}{\aa_r^k}{\aa}{(f_1,\ldots,f_k)}{\Psi_F(f_1,\ldots, f_k)=F(f_1,\ldots, f_k)}$ is continuous,
 \item \label{factiii} for all $f_1,\ldots, f_k\in {\cal A}_r,$ $\|\Psi_F(f_1,\ldots,f_k)\|_*\leq\left(\frac{R}{R-r}\right)^k\|F\|_\infty$.
\end{enumerate}
\end{lemma}
We admit this lemma for the moment and finish the proof of Theorem \ref{jhpd}. We denote by $R$ the radius of convergence of $A_J$, and we apply Lemma \ref{crucial_lemma} to $A_J$. If $f$ belongs to $\Omega$, each of its components belongs to $\aa$ and has norm smaller than $\|f\|$. Therefore, for all $f$ belonging to $\Omega_{\frac R2}:=\{g\in\Omega,\ \|g\|<\frac R2\}$, Lemma \ref{crucial_lemma} implies that each coefficient in the matrix $A_J(f)$ belongs to $\aa$ and have norm smaller that $\|A_J\|_\infty$. Moreover $\left\| \diffp{\overline f}{\overline z}\right\|_*=\sum_{k,l}(1+l)k|f_{k,l}|\leq \|f\|,$
and since $\aa$ is a Banach algebra, we have $\|g\|_*=\left\|A_J(f)\diffp{\overline f}{\overline{z}}\right\|_*\leqs \|A_J\|_\infty \|f\|$, uniformly with respect to $f$ and $J$.\\
Now we notice that $T:\aa^{2n} \to\Omega$ is continuous because for $u=\sum_{k,l}u_{k,l}z^k\overline z^l\in\aa$ we have
$$\int_{[0,\overline z]} u(z,\omega)d\omega =\sum_{k,l}\frac1{l+1} u_{k,l} z^k\overline z^{l+1}$$
and
\begin{align*}
 \sum_{k,l}\frac{k+(l+1)+1+k(l+1)}{l+1} |u_{k,l}|&=\sum_{k,l}\frac{k+1}{l+1}|u_{k,l}|+{(k+1)}|u_{k,l}|\\&\leq 2\|u\|_*.
\end{align*}
Therefore $G=T(g)$ belongs to $\Omega$ and $\|G\|\leq 2\|g\|_*\leqs \|A_J\|_\infty\|f\|$. This proves that $\Phi_J(f)$ belongs to $\Omega$ for all $f\in \Omega_{\frac R2}$ and that $\Phi_J:\Omega_{\frac R2}\to\Omega$ is continuous. Moreover, provided $J$ is close enough to the standard structure, $\|A_J\|_\infty$ is arbitrarily small and so $\|G\|\leq\frac12\|f\|$ which implies that $\Phi_J$ is a small perturbation of the identity. Thus $\Phi_J:\Omega_{\frac R2}\to \Phi_J(\Omega_{\frac R2})$ is continuously invertible.\\
Now Fact (\ref{fact4}) is an immediate consequence of the previous facts because, since $\varphi _2\circ\varphi _1$ is the identity over $\tilde \Omega$, $\varphi_2\circ\Phi_J^{-1}\circ\varphi_1$ is in fact a continuous perturbation of the identity which is arbitrarily small, provided $J$ is close enough to the standard structure.\\
Now, for any sequence $(x_k)_k$ of vectors of $\rr^{2n}$, provided $\|(x_k)_k\|<\frac R2$, $\gamma =\Phi_J^{-1}\circ \varphi _1\circ \left(\varphi _2\circ \Phi_J^{-1}\circ \varphi _1\right)^{-1}\left((x_k)_k\right)$ is a $J$-holomorphic disc such that $\varphi _2(\gamma )=(x_k)_k$. To conclude the proof of the theorem, we have to prove the crucial Lemma \ref{crucial_lemma}.\qed\\[5pt]
{\it Proof of Lemma \ref{crucial_lemma}:} 
The proof of this Lemma is inspired from abstract harmonic analysis and more precisely from Theorem 24D of \cite{Loo}. In order to make the proof clearer, we prove the lemma for $k=2$, that is for $F:\cc^2\to\cc$, analytic over $P(0,R)=\{(z,w)\in\cc^2,\ |z|^2+|w|^2<R^2\}$, $R>0$. The case $k>2$ is a direct generalization of the case $k=2$.

 Let $f=\sum_{k,l} f_{k,l}\zeta^k \overline \zeta^l$ and $g=\sum_{k,l} g_{k,l}\zeta^k \overline \zeta^l$ belonging to $\aa$ be such that $\|f\|_*<r$ and $\|g\|_*<r$ for some $r\in]0,R[$. We write $F(z,w)=\sum_{k,l} F_{k,l}z^kw^l$ and $F(f,g)(\zeta)=\sum_{k,l} h_{k,l}\zeta^k\overline \zeta^l$. In order to determine a suitable expression of $h_{k,l}$ we introduce the Banach algebra
$$\hat\aa=\left\{{(u_{k,l})}_{k,l}\subset \rr^{2n},\ \|(u_{k,l})_{k,l}\|_{\hat*}:= \sum_{k,l}(1+k)|u_{k,l}|<\infty\right\}$$
that we trivially identify to $\aa$ via the application
$$\app\phi\aa{\hat\aa}{\sum_{k,l}u_{k,l}\zeta^k\overline{\zeta}^l}{(u_{k,l})_{k,l}}.$$ We denote by $\mathfrak{e}$ the unit of $\hat\aa$, $\mathfrak e_{k,l}=1$ if $k=l=0$, and $0$ otherwise and we denote by $\indic$ the constant function which equals $1$ so that $\phi(\indic)=\mathfrak{e}$.
We have
\begin{align*}
 h_{k,l}&=\sum_{n,m} F_{n,m}\phi(f^ng^m)_{k,l}
\end{align*}
and using Cauchy's Formula we get:
\begin{align*}
 h_{k,l}&=\sum_{n,m}\frac1{(2i\pi )^2}\int_{\genfrac{}{}{0pt}{}{|\lambda |=R}{|\mu |=R}} \frac{F(\lambda ,\mu )}{\lambda ^{n+1}\mu ^{m+1}} \phi(f^ng^m)_{k,l}d\lambda d\mu.
\end{align*}
Since $\|f\|_*<r$ and $\|g\|_*<r$, we have $\|f^ng^m\|_*<r^{n+m}$ so $|\phi(f^ng^m)_{k,l}|<r^{n+m}$.\\
Therefore if $r<R$, the series $\sum_{n,m}\frac1{\lambda ^n\mu^m}\phi(f^n,g^m)_{k,l}$ converges normally for all $\lambda $ and $\mu $ in $\cc$ such that $|\lambda |=|\mu |=R$, and we can exchange the signs $\int$ and $\sum$ in the last expression of $h_{k,l}$ and which gives
\begin{align*}
 h_{k,l}&=\frac1{(2i\pi )^2}\int_{\genfrac{}{}{0pt}{}{|\lambda |=R}{|\mu |=R}} {F(\lambda ,\mu )}\sum_{n,m}  \left(\phi\left(\frac{f^n}{\lambda ^{n+1}}\right) \phi\left(\frac{g^m}{\mu ^{m+1}}\right)\right)_{k,l}d\lambda d\mu.
\end{align*}
Since $\|f\|_*<r<|\lambda |$ and $\|g\|_*<r<|\mu |$, we have $\left(\lambda \mathfrak e-\phi(f)\right)^{-1}=\sum_{n} \phi\left(\frac{f^n}{\lambda ^{n+1}}\right)$ and 
$\left(\mu \mathfrak e-\phi(g)\right)^{-1}=\sum_{m} \phi\left(\frac{g^m}{\mu ^{m+1}}\right)$
which yields 
 \begin{align*}
  h_{k,l}&=\frac1{(2i\pi )^2}\int_{\genfrac{}{}{0pt}{}{|\lambda |=R}{|\mu |=R}} {F(\lambda ,\mu )}
 \left(\left(\lambda \mathfrak e-\phi(f)\right)^{-1}\left(\mu  \mathfrak e-\phi(g)\right)^{-1}\right)_{k,l}d\lambda d\mu.
\end{align*}
Using this expression of $h_{k,l}$, it follows that
\begin{align*}
 &\sum_{k,l}(k+1)|h_{k,l}|\\
 &\leq \frac{R^2}{(2\pi)^2}\int_0^{2\pi }\hskip-8pt\int_0^{2\pi } |F(Re^{i\theta},Re^{i\varphi })|\sum_{k,l}(k+1)
 \left|\left(\bigl({ Re^{i\theta} \mathfrak e-\phi(f)}\bigr)^{-1}\bigl(Re^{i\varphi } \mathfrak e-\phi(g)\bigr)^{-1}\right)_{k,l}\right|d\theta d\varphi\\
 &\leq  \frac{R^2}{(2\pi)^2}\int_0^{2\pi }\hskip-8pt\int_0^{2\pi } |F(Re^{i\theta},Re^{i\varphi })|
 \left\|\bigl({ Re^{i\theta} \mathfrak e-\phi(f)}\bigr)^{-1}\bigl(Re^{i\varphi } \mathfrak e-\phi(g)\bigr)^{-1}\right\|_{\hat*}d\theta d\varphi.
\end{align*}
Now, using again the identity $\left(\lambda \mathfrak e-\phi(f)\right)^{-1}=\sum_{n} \phi\left(\frac{f^n}{\lambda ^{n+1}}\right)$,  for $\lambda =Re^{i\theta}$ we get $$\left\|\left(\lambda \mathfrak e-\phi(f)\right)^{-1}\right\|_{\hat*}\leq \frac1{R-r}.$$ The same holds true for $g$ and so
$\|F(f,g)\|_*=\sum_{k,l}(k+1)|h_{k,l}|\leq \left(\frac R{R-r}\right)^2\|F\|_\infty$, which implies that $F(f,g)$ belongs to $\aa$ and that (\ref{factiii}) of Lemma \ref{crucial_lemma} holds true.\\
In order to establish the continuity of $\Psi_F$, let $\tilde f$ and $\tilde g$ be two functions of $\aa$ such that $\|\tilde 
f\|_*< r$ and $\|\tilde g \|_*<r$ and set $F(\tilde f,\tilde g)=\sum_{k,l}\tilde h_{k,l}\zeta^k\overline \zeta^l$. We have 
\begin{align*}
 &h_{k,l}-\tilde h_{k,l}=\\
 &\frac1{(2i\pi)^2}
 \int_{\genfrac{}{}{0pt}{}{|\lambda |=R}{|\mu |=R}} {F(\lambda ,\mu )}
 \left(\hskip-3pt\bigl(\lambda \mathfrak e-\phi(f)\bigr)^{-1}\hskip-1pt\bigl(\mu  \mathfrak e-\phi(g)\bigr)^{-1}
 \hskip-2pt-\bigl(\lambda \mathfrak e-\phi(\tilde f)\bigr)^{-1}\hskip-1pt\bigl(\mu  \mathfrak e-\phi(\tilde g)\bigr)^{-1}\hskip-1pt \right)_{k,l}d\lambda d\mu.
\end{align*}
This yields
\begin{align*}
 \|F(f,g)-F(\tilde f,\tilde g)\|_*\leq&\frac{R^2}{4\pi^2}\int_0^{2\pi }\hskip-8pt\int_0^{2\pi } |F(Re^{i\theta},Re^{i\varphi })| \cdot\Bigl\|\bigl({ Re^{i\theta} \mathfrak e-\phi(f)}\bigr)^{-1}\bigl(Re^{i\varphi } \mathfrak e-\phi(g)\bigr)^{-1}\\
& -\bigl({ Re^{i\theta} \mathfrak e-\phi(\tilde f)}\bigr)^{-1}\bigl(Re^{i\varphi } \mathfrak e-\phi(\tilde g)\bigr)^{-1}
 \Bigr\|_{\hat*}d\theta d\varphi.
\end{align*}
We have the following inequality:
\begin{align*}
\Bigl\|\bigl({ Re^{i\theta} \mathfrak e-\phi(f)}\bigr)^{-1}\Bigr.&-\bigl.\bigl(Re^{i\varphi } \mathfrak e-\phi(\tilde f)\bigr)^{-1}\Bigr\|\\
&\leq  \Bigl\|\bigl({ Re^{i\theta} \mathfrak e-\phi(f)}\bigr)^{-1}\cdot \bigl(Re^{i\varphi } \mathfrak e-\phi(\tilde f)\bigr)^{-1}\Bigr\| \cdot\|\phi(f)-\phi(\tilde f)\|_{\hat*}\\
&\leq \frac1{(R-r)^2}\|f-\tilde f\|_*.
\end{align*}
Using this estimate of $\Bigl\|\bigl({ Re^{i\theta} \mathfrak e-\phi(f)}\bigr)^{-1}\Bigr.-\bigl.\bigl(Re^{i\theta} \mathfrak e-\phi(\tilde f)\bigr)^{-1}\Bigr\|$ and the corresponding one for $\Bigl\|\bigl({ Re^{i\varphi } \mathfrak e-\phi(g)}\bigr)^{-1}\Bigr.-\bigl.\bigl(Re^{i\varphi } \mathfrak e-\phi(\tilde g)\bigr)^{-1}\Bigr\|$,
we get
\begin{align*}
 \|F(f,g)-F(\tilde f,\tilde g)\|_*&\leq \frac{R^2}{(R-r)^3}\|F\|_\infty \cdot(\|f-\tilde f\|_*+\|g-\tilde g\|_*)
\end{align*}
which proves that $\Psi_F$ is continuous.\qed

Under a simple growth condition of a sequence $(x_k)_k$, the following corollary gives the existence of a germ of $J$-holomorphic disc $\gamma $ such that for all $k\in\nn$, $\diffp{^k\gamma}{x^k}(0)=x_k$.
\begin{corollary}\mlabel{gjhpd}
 Let $J$ be a real analytic complex structure in a neighborhood of the origin of $\rr^{2n}$, and let $(x_k)_k$ be a sequence of vectors of $\rr^{2n}$. Assume that $x_0$ is close enough to the origin and that there  exists $R>0$ such that for all $k\in\nn^*$, $|x_k|\leq k! R^k$.\\
 Then there exists a germ of  $J$-holomorphic disc $\gamma$ such that for all $k\in\nn$, $\diffp{^k\gamma}{x^k}(0)=x_k$. 
\end{corollary}
\pr Let $r>0$ be the constant given by Theorem \ref{jhpd} and for $\alpha>0$ let $(\tilde x_k)_k$ be the sequence defined by $\tilde x_k=\left(\frac\alpha{R}\right)^k x_k$.\\
Then $\sum_{k=1}^\infty \frac{k+1}{k!}|\tilde x_k|\leq \sum_{k=1}^\infty (k+1)\alpha^k$, so if $|x_0|<r$ and if $\alpha$ is small enough, we have $\sum_{k=0}^\infty \frac{k+1}{k!}|\tilde x_k|<r$. Therefore we can apply Theorem \ref{jhpd} to $(\tilde x_k)_k$ and so there exists a $J$-holomorphic disc $\tilde\gamma:\dd\to\rr^{2n}$ such that  $\diffp{^k\tilde \gamma }{x^k}(0)=\tilde x_k$ for all $k$.\\
Now setting $\gamma(\zeta)=\tilde \gamma( \frac R\alpha \zeta)$ for $\zeta\in\cc$ such that $|\zeta|< \frac \alpha R$, we get a germ of $J$-holomorphic disc such that  $\diffp{^k\gamma}{x^k}(0)=x_k$ for all $k$.\qed

\par\medskip

Given a point $p$ in $M$ and an analytic vector field $X$, we aim at finding a germ of $J$-holomorphic disc $\gamma$ such that $\gamma(0)=p$ and $\diffp{^k\gamma}{x^k}(0)=D_X^{k-1,0}X(p)$ for $k\in\nn^*$ by applying Corollary \ref{gjhpd} to the sequence $(x_k)_k$ defined by $x_0=p$ and $x_k=D_X^{k-1,0}X(p)$, $k\in\nn^*$. For this, we have to prove that the sequence $(D_{X}^{k-1,0}X(p))_k$ satisfies the assumption of this corollary. The estimates needed will be given in the next proposition, but before we state it, we need the following notation.\\
If $X$ is a real analytic vector field in the polydisc ${\cal P}(0,R)=\{(x_1,y_1,\ldots, x_n,y_n)\in\rr^{2n},\ |x_i|<R, |y_i|<R,\ i=1,\ldots, n\}$, we denote by $\tilde X$ the ``polarization'' of $X$. More precisely, if $X(z,\overline z)=\sum_{\nu, \mu} a_{\nu,\mu} z^{\nu}\overline {z}^{\mu}$ where $\nu$ and $\mu$ are multi-indices of $\nn^n$, $z=x+iy$ belongs to $\cc^n$ and $a_{\nu,\mu}$ belongs to $\rr^{2n}$, we set $\tilde X(z,\zeta)=\sum_{\nu, \mu} a_{\nu, \mu} z^{\nu}{\zeta}^{\mu}$ which is defined for all $z$ and $\zeta$ such that $|z_j|,|\zeta_j|< R$, $j=1,\ldots, n$,. We then set
$$c_R(X)=\sup_{|z_j|,|\zeta_j|< R} |\tilde X(z,\zeta)|.$$
We can now state the estimates of $D_X^{k,0}X(0)$:
\begin{proposition}\label{evf}
 Let $X$ be an analytic vector field in a neighborhood of the origin of $\rr^{2n}$. Then for all $R>0$ small enough and for all $k$
$$|D_X^{k,0}X(0)|\leq 2n k! \left(\frac{8nc_R(X)}R\right)^k.$$
\end{proposition}
The proof of Proposition \ref{evf} relies on the following combinatorial Lemmata.
\begin{lemma}\label{lem1}
 For all $\alpha_1,\ldots,\alpha_k\in\nn$ such that $\alpha_1+\ldots+\alpha_k=k$, we have
$$\alpha_1!\ldots\alpha_k!\leq k!.$$
\end{lemma}
\pr Without restriction we assume that $\alpha_1,\ldots,\alpha_r\geq 1$ and $\alpha_{r+1}=\ldots=\alpha_k=0$, $r\leq k$ so that $\alpha_1+\ldots+\alpha_r=k$.\\
The number $k!$ is the number of permutations of the set $\{1,\ldots, k\}$ whereas $\alpha_1!\ldots\alpha_r!$ is the number of permutations of $\{1,\ldots,k\}$ which leave stable each of the sets $\{1,\ldots, \alpha_1\},$ $\{\alpha_1+1,\ldots, \alpha_1+\alpha_2\},\ldots,\{\alpha_1+\ldots+\alpha_{r-1}+1,\ldots, k\}$. Therefore $\alpha_1!\ldots\alpha_k!\leq k!.$\qed
\begin{lemma}\label{lem2}
 There are exactly $\binom{2k-1}{k}$ distinct $k$-tuples $(\alpha_1,\ldots,\alpha_k)\in\nn^k$ such that $\alpha_1+\ldots+\alpha_k=k$.
\end{lemma}
\pr To a $k$-tuple $(\alpha_1,\ldots,\alpha_k)\in\nn^k$ we associate the $k$-tuple $(\beta_1,\ldots,\beta_k)=(\alpha_1+1,\ldots,\alpha_k+1)$ so that $\beta_1+\ldots+\beta_k=2k$. This correspondence is a bijection, and now we count the number of such $k$-tuples $(\beta_1,\ldots,\beta_k)$.\\
The sum $1+1+\ldots+1$ ($2k$-times) equals $2k$ and there are as many ways of writing $\beta_1+\ldots+\beta_k=2k$ as ways of separating the number $1$ in the sum $1+1+\ldots+1=2k$ with $k-1$-sticks, that is $\binom{2k-1}{k-1}$ ways.\qed
\begin{lemma}\label{lem3} For all $k\in\nn$ we have
$$\binom{2k-1}{k-1}\leq \frac{4^k}{\sqrt{\pi k}}.$$
\end{lemma}
\pr 
This is a simple consequence of Stirling's formula.\qed\\
We now prove Proposition \ref{evf}:\\
{\it Proof of Proposition \ref{evf}:} Let $X$ be a real analytic vector field in a neighborhood of the origin in $\rr^{2n}$. Let $R>0$ be such that $X$ is defined and continuous on $\overline{P(0,R)}$.\\
We write $X$ in local coordinates as $X=\sum_{j=1}^{2n} X_j\x j$. In order to establish the proposition, we need an upper bound for the quantity:
\begin{align}
X_{i_1}\x{i_1}\left(X_{i_2}\x{i_2}\left(\ldots\left(X_{i_{k-1}}\diffp{X_{i_k}}{x_{i_{k-1}}}\right)\ldots\right.\right)(0).\label{deriv} 
\end{align}
In the term above, we denote by $j_i$ the cardinal of the set $\{l,\ i_l=i\}$.\\
In (\ref{deriv}), for $i\in\{1,\ldots, 2n\}$ and 
for $\alpha^{(i)}_j$, $j\in\{1,\ldots, k\}$,
 such that $\sum_j\alpha^{(i)}_j=j_i$, the term $\prod_{j=1}^k
 \diffp{^{\alpha_j^{(1)}+\ldots+\alpha_j^{(2n)}}X_{i_j}} {x_1^{\alpha_j^{(1)}}\ldots\partial x_{2n}^{\alpha_{j}^{(2n)}}}(0)$ appears less than $\binom{j_1}{(\alpha_1^{(1)},\ldots,\alpha_k^{(1)})}\ldots\binom{j_{2n}}{(\alpha_1^{({2n})},\ldots,\alpha_k^{({2n})})} $ times.\\
On the other hand, from Cauchy's inequalities, we get
\begin{align*}
 \left|{\diffp{^{\alpha_j^{(1)}+\ldots+\alpha_j^{(2n)}}X_{i_j}} {x_1^{\alpha_j^{(1)}}\ldots\partial x_{2n}^{\alpha_{j}^{(2n)}}}(0)}\right|
&\leq \alpha_j^{(1)}!\ldots\alpha_j^{(2n)}! \frac {c_R(X)} {R^{\alpha_j^{(1)}+\ldots+\alpha_j^{(2n)}}}.
\end{align*}
This yields
\begin{align*}
 \left|X_{i_1}\x{i_1}\left(X_{i_2}\x{i_2}\left(\ldots\left(X_{i_{k-1}}\diffp{X_{i_k}}{x_{i_{k-1}}}\right)\ldots\right.\right)(0)\right|
&\leq\hskip-5pt \sum_{\over{\alpha_1^{(1)}+\ldots+\alpha_k^{(1)}=j_1}{\over{\vdots}{\alpha_1^{(2n)}+\ldots+\alpha_k^{(2n)}=j_{2n}}}}\hskip-5pt
j_1!\ldots j_{2n}! \left(\frac{c_R(X)}{R}\right)^k.
\end{align*}
Applying successively Lemma \ref{lem2}, Lemma \ref{lem3} and Lemma \ref{lem1} we get
\begin{align*}
 \left|X_{i_1}\x{i_1}\left(X_{i_2}\x{i_2}\left(\ldots\left(X_{i_{k-1}}\diffp{X_{i_k}}{x_{i_{k-1}}}\right)\ldots\right.\right)(0)\right|&\leq4^{j_1+\ldots+ j_{2n}}j_1!\ldots j_{2n}! \left(\frac{c_R(X)}{R}\right)^k\\
& \leq k! \left(\frac{4c_R(X)}R\right)^k.
\end{align*}
Now, since there are at most $(2n)^k$ terms $X_{i_1}\x{i_1}\left(X_{i_2}\x{i_2}\left(\ldots\left(X_{i_{k-1}}\diffp{X_{i_k}}{x_{i_{k-1}}}\right)\ldots\right.\right)(0)$ in $D_X^{k,0}X(0)$,
$$|D_X^{k,0}X(0)|\leq 2n k! \left(\frac{8nc_R(X)}R\right)^k$$
which proves the claim.\qed

\section{$J$-holomorphic discs and hypersurfaces}\label{secjhhs}
In this section, given a point $p\in M$ close enough to the origin and a vector $v\in T_p^JM$, we look for a $J$-holomorphic disc $\gamma$ included in $M$ such that $\diffp\gamma  x(p)=v$. 
Our strategy will be the following: we will look for a sequence of vector fields $(X_k)_{k\in\nn^*}$ such that 
\begin{enumerate}[(i)]
 \item For all $k$, the vector field $X_k$ belongs to $T^JM$,
 \item $X_1(p)=v$,
 \item $D_{X_k}^{l,0}X_k(p)=D^{l,0}_{X_1}X_1(p)$ for all $k$ in $\nn^*$ and $l$ in $\nn$,
 \item $X_k$ commutes at order $k$ at $p$.
\end{enumerate}
Then, applying Proposition \ref{evf} to $X_1$ and corollary \ref{gjhpd} to the sequence $(x_k)_k$ defined by $x_0=p$ and $x_k=D^{k-1,0}_{X_1} X_1(p)$ for $k\in\nn^*$, we will get a germ of $J$-holomorphic disc $\gamma$ such that $\gamma(0)=p$ and  $\diffp{^k\gamma}{x^k}(p)=D_{X_1}^{k-1,0} X_1(p)$ for all $k\in\nn^*$, and in particular $\diffp{\gamma }{x}(p)=v$. In order to check that $\gamma $ is included in $M$, we will use  Theorem 1 of \cite{BM}: since $X_k$ commutes at order $k$ at $p$, and since 
$\diffp{^l\gamma}{x^l}(0)=D_{X_k}^{l-1,0} X_k(p)$ for all $l\leq k$, $\gamma$ is tangent to $M$ at $p$ at order $k+1$. Now, since this holds for all $k$, this implies that $\gamma $ is tangent to $M$ at $p$ at any order. As $M$ and $\gamma $ are real analytic, $\gamma$ is in fact included in $M$. In order to construct such a sequence of vector fields, we will need the following propositions.
\begin{proposition}\label{prop1}
 Let $X$ be any real analytic vector field in a neighborhood of the origin. Then the following assertions are equivalent~:
 \begin{enumerate}[(i)]
  \item \label{i} For all $0\leq k\leq k_0-2$ and  all $X_1,\ldots, X_k\in \{X,JX\}$, $X_1\cdot X_2\ldots X_k\cdot[X,JX](0)=0$,
  \item \label{ii} $X$ commutes at order $k_0$ at $0$. 
 \end{enumerate}
Moreover, in that case, for all integers $k$ and $l$ such that $0\leq l\leq k\leq k_0-1$ the following holds:
\begin{align}
X_1\ldots X_l\cdot\left[X_{l+1},\left[\ldots\left[X_k,[X,JX]]\right.\ldots\right.\right](0)&=\left[X_1,\left[X_2,\ldots\left[X_k,[X,JX]\right]\ldots\right.\right](0). \mlabel{eq10}
\end{align}
\end{proposition}
\begin{remark} \rm
 We point out that (\ref{eq10}) holds true for $0\leq  k\leq k_0-1$ even if (\ref{i}) is true only for $0\leq k\leq k_0-2$.
\end{remark}
\noindent{\it Proof of proposition \ref{prop1}:} The equivalence of (\ref{i}) and (\ref{ii}) was shown in \cite{BM}, proposition 13. The second assertion is proved by induction on $l$, the case $l=0$ being trivial.\\
We assume that the identity holds true for $l\leq k\leq k_0-2$ and we denote by $Y$ the vector field $Y=[X_{l+2},[\ldots,[X_k,[X,JX]]\ldots]]$. Then, since $\tilde X_1\ldots \tilde X_l \cdot Y(0)=0$ for all $\tilde X_1,\ldots \tilde X_l\in\{X,JX\}$, we have
\begin{align*}
 &X_1\ldots X_l[X_{l+1},[X_{l+2},\ldots[X_k,[X,JX]]\ldots](0)\\
 &= X_1\ldots X_{l}\cdot X_{l+1}\cdot Y(0)-X_1\ldots X_l \cdot Y\cdot X_{l+1}(0)\\
 &=X_1\ldots X_{l}\cdot X_{l+1}\cdot Y(0).
\end{align*}
\qed
\begin{proposition}\label{prop2}
 Let $X$ be a vector field commuting at order $k$ at $0$, $k\geq 2$. 
 Then for all $X_1,\ldots, X_{k}\in\{X,JX\}$, we have $$[X_1,[\ldots,[X_{k},[X,JX]]]](0)=\underbrace{[JX,[\ldots,[JX}_q,\underbrace{[X,[\ldots,[X}_p,[X,JX]\ldots](0)$$
 where $p=\#\{i,\ X_i=X\}$ and $q=\#\{i,\ X_i=JX\}.$ 
\end{proposition}
\pr It suffices to show that for all $l\leq k-2$ we have
\begin{align*}
 &[X_1,[\ldots,[X_l,[X,[JX,[X_{l+3},[\ldots,[X,JX]\ldots](0)\\
 &\hskip50pt = [X_1,[\ldots,[X_l,[JX,[X,[X_{l+3},[\ldots,[X,JX]\ldots](0).
\end{align*}
We set $Y=[X_{l+3},[\ldots,[X,JX]\ldots]$. From Jacobi identity we have:
\begin{align}
[X,[JX,Y]]&=[JX,[X,Y]]+[[X,JX],Y]. \mlabel{eq12}
\end{align}
From Proposition \ref{prop1}, for all $0\leq \tilde l\leq l$ and all  
$\tilde X_1,\ldots, \tilde X_{\tilde l}\in\{X,JX\}$, we have  $\tilde X_1\ldots\tilde X_{\tilde l}\cdot [X,JX](0)=0$ which yields
\begin{align}
 [X_1,\ldots,[X_l,[X,JX]\cdot Y]\ldots](0)&=0\label{eq11}
\end{align}
Analogously, it follows from Proposition \ref{prop1} that
$[X_1,\ldots,[X_l,Y\cdot[X,JX]]\ldots](0)=0$ which with (\ref{eq11}) gives
\begin{align*}
 [X_1,\ldots,[X_l,[[X,JX],Y]]\ldots](0)=0
\end{align*}
and with (\ref{eq12}) we end the proof of the proposition.\qed

Taking our strategy into account, the proof of Theorem \ref{alg-lie} will be a direct consequence of the following lemma:
\begin{lemma}\label{pfreeman}
Let $X\in T^JM$ be a real analytic vector field in a neighborhood of $0\in M$ such that  $\ll(X)$, the Lie algebra generated by $X$, is included in $T^JM$.\\
Then there exists a sequence $(X_k)_{k\in\nn^*}$ of vector fields in $T^JM$ such that 
\begin{enumerate}[(a)]
 \item\label{a} $D_{X_k}^{l,0}X_k(0)=D^{l,0}_{X}X(0)$ for all $k$ in $\nn^*$ and $l$ in $\nn$,
 \item \label{b} for all $k\in\nn^*$,  $X_k$ commutes at order $k$ at $0$.
\end{enumerate}
\end{lemma}
\pr We set $X_1=X$ and by induction on $k$ we construct a sequence $(X_k)_k$ of vector fields which satisfy (\ref{a}) and (\ref{b}) and which are such that for all $k$, $X_k$ belongs to $\ll(X)$.  Assuming $X_k$ constructed, we set 
$$Y_{p,q}=\underbrace{[JX_k,[\ldots,[JX_k,}_{q}\underbrace{[X_k,[\ldots,[X_k}_{p}[X_k,JX_k]\ldots ]$$
for all $p$ and $q$ such that $p+q=k-1$, and we look for $X_{k+1}$ as
\begin{align*}
 X_{k+1}&= X_k+\sum_{p+q=k-1} (a_{p,q}+b_{p,q}J )Y_{p,q},
\end{align*}
where $a_{p,q}$ and $b_{p,q}$ are real analytic functions which vanish at order $k-1$ at the origin. Hence $X_{k+1}$ already belongs to $\ll(X)$ and we have to show that $a_{p,q}$ and $b_{p,q}$ can be chosen so that (\ref{a}) and (\ref{b}) hold. 
We compute $[X_{k+1},JX_{k+1}]$:
\begin{align*}
 &[X_{k+1},JX_{k+1}]=\\
 &=[X_k,JX_k]+\sum_{p+q=k-1} \Bigl(\bigl(X_k(a_{p,q})-JX_k(b_{p,q})\bigr) JY_{p,q}+\bigl(-X_k(b_{p,q})-JX_k(a_{p,q})\bigr) Y_{p,q}\Bigr)\\
 &+\sum_{p+q=k-1} (a_{p,q}^2+b_{p,q}^2) [Y_{p,q},JY_{p,q}]\\
 &+\sum_{p+q=k-1} \bigl(a_{p,q} \left([X_k,JY_{p,q}]+[Y_{p,q},JX_k]\right)-b_{p,q}\left([X_k,Y_{p,q}]+[JX_k,JY_{p,q}]\right)\bigr)\\
 &+\sum_{\over{p+q=k-1}{r+s=k-1}}\bigl(a_{r,s}Y_{r,s}(a_{p,q})+b_{r,s}JY_{r,s}(a_{p,q})-a_{r,s}JY_{r,s}(b_{p,q})+b_{r,s}Y_{r,s}(b_{p,q})\bigr)JY_{p,q}\\	
 &+\sum_{\over{p+q=k-1}{r+s=k-1}}\bigl(-a_{r,s}Y_{r,s}(b_{p,q})-b_{r,s}JY_{r,s}(b_{p,q})-a_{r,s}JY_{r,s}(a_{p,q})+b_{r,s}Y_{r,s}(a_{p,q})\bigr)Y_{p,q}.
\end{align*}
Since $a_{p,q}$ and $b_{p,q}$ vanish at order $k-1$ at the origin, and since $X_k$ commutes at order $k$ at $0$,  for all $X^{(1)},\ldots, X^{(k-2)}\in\{X_{k+1},JX_{k+1}\}$ we have
 $$X^{(1)}\ldots X^{(k-2)}[X_{k+1}, JX_{k+1}](0)=0.$$
 It then follows from Proposition \ref{prop1} that $X_{k+1}$ already commutes at order $k$ at $0$.  Moreover, $D_{X_{k+1}}^{r,0} X_{k+1}(0)=D_{X_{k}}^{r,0} X_{k}(0)$ for all $r\leq k-1$ and for all choices of functions $a_{p,q}$ and $b_{p,q}$, vanishing at order $k-1$ at $0$.\\
We now choose these functions so that $X_{k+1}$ commutes at order $k+1$ at the origin and so that
$D_{X_{k+1}}^{r,0} X_{k+1}(0)=D_{X_{k}}^{r,0} X_{k}(0)$ for all $r\geq k$.

It follows from Proposition \ref{prop1} that $Y_{r,s}(0)=D^{r,s}_{X_k} [X_k,JX_k](0)$ for all $r$ and $s$ such that $r+s=k-1$. Since the functions $a_{p,q}$ and $b_{p,q}$ vanish at order $k-1$ at the origin, we have for all $r$ and $s$ such that $r+s=k-1$
\begin{align*}
 D^{r,s}_{X_{k+1}}[X_{k+1},JX_{k+1}](0)=
 &Y_{r,s}(0)+\sum_{p+q=k-1} \left(\left(D^{r+1,s}_{X_k} a_{p,q}(0)-D^{r,s+1}_{X_k} b_{p,q}(0)\right) JY_{p,q}(0)\right)\\
 &+\sum_{p+q=k-1} \left(\left(-D^{r+1,s}_{X_k}b_{p,q}(0)-D^{r,s+1}_{X_k} a_{p,q}(0)\right) Y_{p,q}(0)\right).
\end{align*}
Perhaps after a change of coordinates if needed, we can assume that $X_k=\x1$ and that $J(0)$ is still the standard structure. We therefore have $JX_k=\sum_{k=1}^n J_{2k-1,1} \x k+ J_{2k} \y k$ and since $J(0)$ is the standard structure, $JX_k(0)=\y 1$. Now, since $a_{p,q}$ and $b_{p,q}$ vanishes at $0$ at order $k-1$, neither $X_k$ nor $JX_k$ is differentiated in $D^{r,s}_{X_{k}}a_{p,q}(0)$ and $D^{r,s}_{X_{k}}b_{p,q}(0)$ and we therefore have
\begin{align}
 \label{eq20} &D^{r,s}_{X_{k+1}}[X_{k+1},JX_{k+1}](0)=\\
 \nonumber&=Y_{r,s}(0)+\sum_{p+q=k-1} \left(\left(\diffp{^{k}a_{p,q}}{x_1^{r+1}\partial y_1^s}(0)-\diffp{^{k}b_{p,q}}{x_1^{r}\partial y_1^{s+1}}(0)\right) JY_{p,q}(0)\right)\\
 \nonumber &+\sum_{p+q=k-1} \left(\left(-\diffp{^{k}b_{p,q}}{x_1^{r+1}\partial y_1^s}(0)-\diffp{^{k}a_{p,q}}{x_1^{r}\partial y_1^{s+1}}(0)\right) Y_{p,q}(0)\right).
\end{align}
We search for $a_{p,q}$ and $b_{p,q}$ as homogeneous polynomials in $x_1$ and $y_1$ of degree $k$ without the term $x_1^{k}$ such that
$$\forall (r,s)\in\nn,\ r+s=k-1, \ \begin{cases}
\diffp{^{k}a_{p,q}}{x_1^{r+1}\partial y_1^s}(0)-\diffp{^{k}b_{p,q}}{x_1^{r}\partial y_1^{s+1}}(0)&=0,\\
\diffp{^{k}b_{r,s}}{x_1^{r+1}\partial y_1^s}(0)+\diffp{^{k}a_{p,q}}{x_1^{r}\partial y_1^{s+1}}(0)&=0\text{ if } (r,s)\neq (p,q),\\
\diffp{^{k}b_{p,q}}{x_1^{r+1}\partial y_1^s}(0)+\diffp{^{k}a_{p,q}}{x_1^{r}\partial y_1^{s+1}}(0)&=1.
\end{cases}$$
With such a choice, (\ref{eq20}) gives $D^{r,s}_{X_{k+1}}[X_{k+1},JX_{k+1}]=0$, for all $r$ and $s$ such that $r+s=k-1$, which together with Propositions \ref{prop1} and \ref{prop2} proves that $X_{k+1}$ commutes at order $k+1$ at 0. So (\ref{b}) holds true.\\
Now, we notice that $D_{X_k}^{r,0} a_{p,q}(0)=D_{X_k}^{r,0} b_{p,q}(0)=0$ for all $r\in\nn$, and all $(p,q)$. Therefore $D_{X_{k+1}}^{r,0}X_{k+1}(0)=D_{X_{k}}^{r,0}X_{k}(0)$ for all $r\in\nn$, and so (\ref{a}) also holds true.\qed
\par\medskip
We deduce from Theorem \ref{alg-lie} the following generalization of Freeman's Theorem:
\begin{theorem}\label{freeman}
 If $J$ is integrable, then for all $X\in\ker \cll$ and all $p$ in a neighborhood of $0$, there exists a $J$-holomorphic disc $\gamma$ such that $\gamma(0)=p$ and $\diffp{\gamma }{x}(0)=X(p)$.
\end{theorem}
\pr Let $X$ be a vector field in $\ker \cll$. We show that the Lie algebra generated by $X$ and $JX$ is included in $T^JM$ by proving by induction that all Lie brackets we can form with $X$ and $JX$ belong to $\ker\cll$.

In the integrable case, we have the following characterization of $\ker \cll$ (see \cite{Fre1}):
\begin{align}
\label{eq21} \ker \cll&=\{X\in T^JM,\ \forall Y\in T^JM,\ [X,Y]\in T^JM\}. 
\end{align}
Assume that all Lie brackets of length at most $k$ belongs to $\ker\cll$. Let $X_p$ and $X_q$ be two Lie brackets of length $p$ and $q$ respectively, with $p$ and $q$ such that $p+q=k+1$. Using (\ref{eq21}), we show that $[X_p,X_q]$ and $[X_p,JX_q]$ belong to $\ker\cll$.\\
Jacobi's identity gives for all $Y\in T^JM$:
\begin{align*}
 [[X_p,JX_q],Y]&=[X_p,[JX_q,Y]]+[JX_q,[Y,X_p]].
\end{align*}
The vector field $JX_q$ belong to $\ker\cll$ so (\ref{eq21}) implies that $[JX_q,Y]$ belongs $T^JM$ and so, again (\ref{eq21}) implies that $[X_p,[JX_q,Y]]$ belongs to $T^JM$. Analogously, $[JX_q,[Y,X_p]]$ belongs to $T^JM$ and so $[[X_p,JX_q],Y]$ is a element of $T^JM$ for all $Y$ in $T^JM$, which with (\ref{eq21}) implies that $[X_p,J X_q]$ belongs to $\ker \cll$.\\
Analogously, $[X_p, X_q]$ belongs to $\ker\cll$. Therefore all Lie brackets of $X$ and $JX$ belong to $\ker\cll$ and so to $T^JM$. Applying Theorem \ref{alg-lie} ends the proof of Theorem \ref{freeman}.\qed
\par\medskip
Another example of situation where Theorem \ref{alg-lie} can be useful is the following example:
\begin{example}\label{ex1}
\rm Let $\varphi:\rr^8\to\rr$ be the map defined by $\varphi(x_1,y_1,x_2,y_2,x_3,y_3,x_4,y_4)=y_1$ and let $M$ be the set $M=\{z\in\rr^8,\ \varphi(z)=0\}$. We  define the eight following vector fields 
\begin{align*}
 L_1&=\diffp{}{x_1},			&L_2&=\y{1},\\
L_3&=\x2-\frac12y^2_3\diffp{}{x_1},	&L_4&=\diffp{}{y_2}+(-2y_3x_3+x_2)\diffp{}{x_1},\\
L_5&=\x{3}-y_2y_3\x{1},			&L_6&=\y3+x_3\x2-\left(\frac{x_3y_3^2}2+x_3y_2\right)\x1,\\
L_7&=\x4+y_4\y4				&L_8&=\y4.
\end{align*}
and the complex structure $J$ they induce by setting 
\begin{align*}
 JL_1&=L_2,&JL_3&=L_4,&JL_5&=L_6, &JL_7=L_8.
\end{align*}
This example is derived from the example of Section \ref{contre-exemple}. Again $J(0)$ is the standard complex structure $J_0$ and the tangent space $TM$ is spanned over $\rr$ by $L_1,\ L_3,\ L_4,\ L_5,\ L_6,\ L_7$ and $L_8$ and $T^JM$ is spanned over $\cc$ by $L_3$, $L_5$ and $L_7$.\\
When we compute the Lie brackets of the $L_i$'s in the complex tangent bundle, we get:
\begin{align*}
 [L_3,L_4]&=L_1,	\\
 [L_3,L_5]&=0,			& [L_4,L_5]&=y_3\x1,		\\
 [L_3,L_6]&=y_3\x1,		& [L_4,L_6]&=0,			&[L_5,L_6]&=L_3,\\ 
 [L_3,L_7]&=0, 			& [L_4,L_7]&=0,			&[L_5,L_7]&=0,			&[L_6,L_7]&=0,\\
 [L_3,L_8]&=0		,	& [L_4,L_8]&=0,			&[L_5,L_8]&=0	,		&[L_6,L_8]&=0,	&[L_7,L_8]&=-L_8.
\end{align*}
Therefore 
\begin{align*}
 &\hat \cll_\varphi(L_3,L_5)=0,& & \hat \cll_\varphi(L_3,L_7)=0,& &\hat \cll_\varphi(L_5,L_7)=0,\\
 &\hat\cll_\varphi(L_3,JL_3)= d\varphi(L_2)=1,& & \hat \cll_\varphi(L_5,JL_5)=d\varphi(L_3)=0,& &\hat \cll_\varphi(L_7,JL_7)=d\varphi(-L_8)=0,
 \end{align*}
and so, $M$ is pseudoconvex and $\ker \cll$ is spanned by $L_5$ and $L_7$.

As in Section \ref{contre-exemple}, there is no $J$-holomorphic disc $\gamma$ included in $M$ such that $\gamma(0)=0$ and $\diffp\gamma x(0)=\x3$ but  Theorem \ref{alg-lie} yields the existence of a $J$-holomorphic disc $\gamma$ included in $M$, such that $\gamma(0)=0$ and $\diffp \gamma x(0)=L_7(0)$.
\end{example}
\par\medskip
We now prove Theorem \ref{subbundle} by applying our strategy with the following lemma:
\begin{lemma}\label{ppfreeman}
Let $\mathbb{L}$ be a subbundle of $T^JM$ such that if $X\in \ll$ commutes at order $k$ at $0$, then for all $X_1,\ldots, X_{k+1}\in \{X,JX\}$, $[X_1,[\ldots,[X_{k},X_{k+1}]\ldots]](0)$ belongs to $\ll_0$.\\
Then for all $X\in\ll$, there exists a sequence $(X_k)_k$ of vector fields in $\ll$ such that  
\begin{trivlist}
 \item[$(\alpha)$]\label{alpha} $D_{X_k}^{l,0}X_k(0)=D^{l,0}_{X}X(0)$ for all $k$ in $\nn^*$ and $l$ in $\nn$,
 \item[$(\beta)$] \label{beta} $X_k$ commutes at order $k$ at $0$.
\end{trivlist}
\end{lemma}
\pr We proceed by induction on $k$. We set $X_1=X$ and we assume $X_k$ constructed. Maybe after a change of coordinates if needed, we can assume that $X_k=\x1$ and $JX_k(0)=\y1$. We set
\begin{align*}
 X_{k+1}&= X_k+\sum_{p=1}^\ell (a_{p}+Jb_{p} )L_{p}
\end{align*}
where $L_1,\ldots, L_\ell$ is a basis of $\ll$ and where $a_{p}$ and $b_{p}$, $p=1,\ldots, \ell$,  are homogeneous polynomials of degree $k$ in the variables $x_1$ and $y_1$, without the term $x^{k}_1$. Therefore $X_{k+1}$ belongs to $\ll$,  
$D_{X_{k+1}}^{l,0}X_{k+1}(0)=D^{l,0}_{X_1}X_1(0)$ for all $l$ and $X_{k+1}$ commutes at order $k$ at $0$.\\
Analogously to the proof of Lemma \ref{pfreeman}, we now choose $a_{p}$ and $b_{p}$ so that $X_{k+1}$ commutes at order $k+1$ at $0$. We compute
\begin{align*}
 [X_{k+1},JX_{k+1}]&=[X_k,JX_k]+\sum_{p=1}^\ell(a_{p}^2+b_{p}^2) [L_{p},JL_{p}]\\
 &+\sum_{p=1}^\ell \Bigl(\bigl(X_k(a_{p})-JX_k(b_{p})\bigr) JL_p+\bigl(-X_k(b_{p})-JX_k(a_{p})\bigr) L_{p}\Bigr)\\
 &+\sum_{p=1}^\ell \bigl(a_{p} \left([X_k,JL_p]+[L_{p},JX_k]\right)-b_{p}\left([X_k,L_{p}]+[JX_k,JL_{p}]\right)\bigr)\\
 &+\sum_{p,q=1}^{\ell}\bigl(a_{q}L_{q}(a_{p})+b_{q}JL_{q}(a_{p})-a_{q}JL_{q}(b_{p})+b_{q}L_{q}(b_{p})\bigr)JL_{p}\\	
 &+\sum_{p,q=1}^{\ell}\bigl(-a_{q}L_{q}(b_{p})-b_{q}JL_{q}(b_{p})-a_{q}JL_{q}(a_{p})+b_{q}L_{q}(a_{p})\bigr)L_{p}.
\end{align*}
On the one hand, 
the functions $a_{p}$ and $b_{p}$ vanish at order $k-1$ at the origin. Therefore, for all $r$ and $s$ such that $r+s=k-1$
\begin{align*}
 D^{r,s}_{X_{k+1}}[X_{k+1},JX_{k+1}](0)&= D^{r,s}_{X_{k}}[X_{k},JX_{k}](0)+\sum_{p=1}^\ell \Bigl(\bigl(D^{r+1,s}_{X_k} a_{p}(0)-D^{r,s+1}_{X_k} b_{p}(0)\bigr) JL_{p}(0)\Bigr)\\
 &+\sum_{p=1}^\ell \Bigl(\bigl(-D^{r+1,s}_{X_k}b_{p}(0)-D^{r,s+1}_{X_k} a_{p}(0)\bigr) L_{p}(0)\Bigr).\\
\end{align*}

On the other hand,  Proposition \ref{prop1}  implies that for all $r$ and $s$ such that $r+s=k-1$  
$$D^{r,s}_{X_k}[X_k, JX_k](0)=[\underbrace{JX_k,[\ldots, [JX_k}_s [\underbrace{X_k,\ldots[X_k}_r,[X_k, JX_k]]\ldots](0)$$
and so $D^{r,s}_{X_k}[X_k, JX_k](0)$ belongs to $\ll_0$. Hence there exists $\alpha^{r,s}_1,\ldots, \alpha^{r,s}_\ell\in\cc$ such that
$$D^{r,s}_{X_k}[X_k, JX_k](0)=\sum_{p=1}^\ell \alpha^{r,s}_p L_p(0)$$
 and so 
 \begin{align*}
 D^{r,s}_{X_{k+1}}[X_{k+1},JX_{k+1}](0)&=\sum_{p=1}^\ell \alpha^{r,s}_p L_p(0)+\sum_{p=1}^\ell \Bigl(\bigl(D^{r+1,s}_{X_k} a_{p}(0)-D^{r,s+1}_{X_k} b_{p}(0)\bigr) JL_{p}(0)\Bigr)\\
 &+\sum_{p=1}^\ell \Bigl(\bigl(-D^{r+1,s}_{X_k}b_{p}(0)-D^{r,s+1}_{X_k} a_{p}(0)\bigr) L_{p}(0)\Bigr).\\
\end{align*}
Now we choose the functions $a_{p}$ and $b_{p}$ as in Lemma \ref{pfreeman} such that
$$\forall (r,s)\in\nn,\ r+s=k-1, \ \begin{cases}
\diffp{^{k}a_{p}}{x_1^{r+1}\partial y_1^s}(0)-\diffp{^{k}b_{p}}{x_1^{r}\partial y_1^{s+1}}(0)&=-\im(\alpha^{r,s}_p)\\
\diffp{^{k}b_{r,s}}{x_1^{r+1}\partial y_1^s}(0)+\diffp{^{k}a_{p,q}}{x_1^{r}\partial y_1^{s+1}}(0)&=\re(\alpha^{r,s}_p)
\end{cases}.$$
After solving this system of $2k$ equations and $2k$ variables, we get a vector field $X_{k+1}$ such that $D_{X_{k+1}}^{r,s}[X_{k+1}, JX_{k+1}](0)=0$ for all $r$ and $s$ such that $r+s=k-1$. Since $X_{k+1}$ already commutes at order $k$, Propositions \ref{prop1} and \ref{prop2} then imply that in fact $X_{k+1}$ commutes at order $k+1$ at $0$.\qed

As a corollary of Theorem \ref{subbundle}, we get the following Theorem already proved by Kruzhilin and Sukhov in \cite{KS} in the smooth case:
\begin{theorem}\label{KS}
 If $\ker\cll =T^JM$, then for all $p\in M$ and all $v\in T^JM$, there exists a germ of $J$-holomorphic disc $\gamma $ such that $\gamma (0)=p$, $\diffp{\gamma }{x}(p)=v$ and $\gamma(\dd)\subset M$.
\end{theorem}
\pr We prove that if $X\in T^JM$ commutes at order $k$ at $0$, then for all $X_1,\ldots, X_{k-1}\in \{X,JX\}$, $[X_1,[\ldots,[X_{k-1},[X,JX]]\ldots](0)$ belongs to $T_0^JM$ and we apply Theorem \ref{subbundle}.\\
On the one hand, since $M$ is Levi flat, $[X,JX]$ belongs to $T^JM$ and $d\varphi (J[X,JX])$ vanishes identically.\\
On the other hand, since $X$ commutes at order $k$ at $0$ 
\begin{align*}
 X_1\cdot X_2\ldots X_{k-1}\cdot d\varphi (J[X,JX])(0)
 &=d\varphi (J(X_1\cdot X_2\ldots X_{k-1}\cdot [X,JX]))(0)\\
 &= d\varphi (J[X_1,[X_2,\ldots [X_{k-1},[X,JX]\ldots])(0)\\
\end{align*}
and so $[X_1,\ldots [X_{k-1},[X,JX]\ldots](0)$ belongs to $JT_0M$. Considering $d\varphi([X,JX])$, one can show analogously that $[X_1,\ldots [X_{k-1},[X,JX]\ldots](0)$ belongs to $T_0M$ and so finally to $T^J_0M.$\qed
\par\medskip

Let us notice that the above theorem can in fact be generalized in the following sense. If $\ll$ is a complex bundle included in $\ker \cll$ such that
\begin{itemize}
 \item for all $X\in\ll$, $[X,JX]$ belongs to $\ll$,
 \item there exist $\varphi_1,\ldots, \varphi_k$ such that for all $X\in T^JM$, $X\cdot\varphi _1=\ldots=X\cdot \varphi _k=0$ if and only if $X$ belongs to $\ll$,
\end{itemize}
then the previous proof shows that for any $X$ in $\ll$, there exists   a $J$-holomorphic disc $\gamma $ such that $\gamma (0)=p$, $\diffp{\gamma }{x}(p)=v$ and $\gamma(\dd)\subset M$. However, as we will see in the next example, Theorem \ref{subbundle} is more general.
\begin{example}\label{ex2}\rm
 Let $\varphi:\rr^8\to\rr$ be the map defined by $\varphi(x_1,y_1,x_2,y_2,x_3,y_3,x_4,y_4)=y_1$ and let $M$ be the set $M=\{z\in\rr^8, \ \varphi(z)=0\}$. We also define the eight following vector fields 
\begin{align*}
 L_1&=\diffp{}{x_1},	&L_2&=\y{1},\\
L_3&=\x2-\frac12y^2_3\diffp{}{x_1},&L_4&=\diffp{}{y_2}+(-2y_3x_3+x_2)\diffp{}{x_1},\\
L_5&=\x{3}-y_2y_3\x{1},		&L_6&=\y3+x_3\x2-\left(\frac{x_3y_3^2}2+x_3y_2\right)\x1,\\
L_7&=\x4+y_4\y4+y_1y_4 L_5	&L_8&=\y4+y_1y_4L_6.
\end{align*}
and the complex structure $J$ they induce by setting 
\begin{align*}
 JL_1&=L_2,&JL_3&=L_4,&JL_5&=L_6, &JL_7=L_8.
\end{align*}
Again this example is derived from the example of Section \ref{contre-exemple}. The tangent space $TM$ is spanned over $\rr$ by $L_1,\ L_3,\ L_4,\ L_5,\ L_6,\ L_7$ and $L_8$ and $T^JM$ is spanned over $\cc$ by $L_3$, $L_5$ and $L_7$.\\
When we compute the Lie brackets of the $L_i$'s in the complex tangent bundle, we get:
\begin{align*}
 [L_3,L_4]&=L_1,	\\
 [L_3,L_5]&=0,			& [L_4,L_5]&=y_3\x1,		\\
 [L_3,L_6]&=y_3\x1,		& [L_4,L_6]&=0,			&[L_5,L_6]&=L_3,\\ 
 [L_3,L_7]&=0, 			& [L_4,L_7]&=y_1y_3y_4\x1,	&[L_5,L_7]&=0,			&[L_6,L_7]&=y_1y_4L_3,\\
 [L_3,L_8]&=y_1y_3y_4\x1,	& [L_4,L_8]&=0,			&[L_5,L_8]&=y_1y_4L_3,		&[L_6,L_8]&=0,	\\
\end{align*}
and
$$[L_7,L_8]=-y_1L_5+y_1^2y_4^2L_3+2y_1y_4L_6-L_8.$$
Therefore
\begin{align*}
 \hat \cll_\varphi(L_3,L_5)&=0,\\
 \hat \cll_\varphi(L_3,L_7)&=0,\\
\hat \cll_\varphi(L_5,L_7)&=0,\\
 \hat\cll_\varphi(L_3,JL_3)&= d\varphi(L_2)=1,\\
 \hat \cll_\varphi(L_5,JL_5)&=d\varphi(L_3)=0,\\
\hat \cll_\varphi(L_7,JL_7)&=d\varphi(y_1^2y_4^2L_4-y_1y_4L_5-y_1L_6-L_8)=0,\\
 \end{align*}
and so, $M$ is pseudoconvex and $\ker \cll$ is spanned by $L_5$ and $L_7$. As in Section \ref{contre-exemple}, there is no $J$-holomorphic disc $\gamma$ included in $M$ such that $\gamma(0)=0$ and $\diffp\gamma x(0)=\x3$. However, Theorem \ref{subbundle} can be applied to $L_7$ with $\ll$ the subbundle generated by $L_7$. In order to do this, it suffices to check that for all $X_1,\ldots, X_k\in \{L_7,L_8\}$, the Lie bracket $[X_1,[\ldots,[X_k,[L_7,L_8]\ldots](0)$ is a linear combination of $L_7$ and $L_8$. By induction on $k$, we get
$$[X_1,[\ldots,[X_k,[L_7,L_8]\ldots]=y_1L+\alpha L_8$$
where $L$ is a vector field and $\alpha$ belongs to $\{-1,0,1\}$, both depending on the sequence $X_1,\ldots, X_k$. Therefore $[X_1,[\ldots,[X_k,[L_7,L_8]\ldots](0)$ belongs to the subbundle generated by $L_7$ at $0$ and Theorem \ref{subbundle} yields the existence of a $J$-holomorphic disc $\gamma$ included in $M$, such that $\gamma(0)=0$ and $\diffp \gamma x(0)=L_7(0)$. 
\end{example}


\begin{thebibliography}{999}
\bibitem{BM} J.-F. Barraud, E. Mazzilli: \it Regular type of real hyper-surface in (almost) complex manifolds,\rm\ Math. Z. 248, 757--772 (2004)
\bibitem{Fre1} M. Freeman: {\it Local Complex foliation of Real Submanifolds,} Math. Ann. 209, 1-30 (1974).
\bibitem{IR} S. Ivashkovich, J.-P. Rosay: {\it  Schwarz-type lemmas for solutions of $\overline\partial$-inequalities and complete hyperbolicity of almost complex manifolds,} Ann. Inst. Fourier (Grenoble) 54 (2004), no. 7, 2387-2435 (2005).
\bibitem{IS} S. Ivashkovich, A. Sukhov: \it Schwarz reflection principle, boundary regularity and compactness for $J$-complex curves,\rm
Ann. Inst. Fourier 60, no. 4, 1489--1513 (2010).
\bibitem{KS} N. Kruzhilin, A. Sukhov: {\it Pseudoholomorphic discs attached to CR-submanifolds of almost complex spaces,} Bull. Sci. Math. 129, no. 5, 398-414 (2005).
\bibitem{Loo} L. H. Loomis: {\it  An introduction to abstract harmonic analysis,} \ D. Van Nostrand Company, Inc., Toronto-New York-London, (1953).
\bibitem{Sik} J.-C. Sikorav: {\it Some properties of holomorphic cueves in almost complex manifolds,} Prog. Math. 117, 165-189, Birkhauser, Basel (1994).
\end{thebibliography}
\end{document}